%% file: HistoryoftheShift.tex
\documentclass{amsart}
\usepackage{amssymb}
\usepackage{amsmath}
\usepackage{graphicx}

\newcommand{\ie}{\textit{i.e.,}}
\newcommand{\apriori}{\textit{a priori}}
\newcommand{\alp}{\ensuremath{\alpha}}
\newcommand{\eps}{\ensuremath{\varepsilon}}
\newcommand{\omeg}{\ensuremath{\omega}}
\newcommand{\epsgo}{\ensuremath{\eps>0}}
\newcommand{\as}[1]{\ensuremath{a_{#1}}}
\newcommand{\bs}[1]{\ensuremath{b_{#1}}}
\newcommand{\ns}[1]{\ensuremath{n_{#1}}}
\newcommand{\recip}[1]{\ensuremath{\frac{1}{#1}}}
\newcommand{\fof}[1]{\ensuremath{f(#1)}}
\newcommand{\ypp}{\ensuremath{y^{\prime\prime}}}

\begin{document}

\title[Genesis of Symbolic Dynamics]{On the Genesis of Symbolic Dynamics\\as We Know It}
\author{Ethan M. Coven}
	\address{Department of Mathematics, Wesleyan University, Middletown CT}
        \email{ecoven@wesleyan.edu}
\author{Zbigniew H. Nitecki}
	\address{Department of Mathematics, 
        Tufts University, Medford, MA 02155, USA}
        \email{zbigniew.nitecki@tufts.edu}
\keywords{Symbolic dynamics, shift automorphism, recurrence, transitivity, minimal set, geodesics}
\subjclass[2000]{Primary 37B10; Secondary 01A60}
\input{abstract}

\maketitle
\input{ShiftText}

\bibliographystyle{amsalpha}
\bibliography{DynamicsHistory}
 
 \end{document}

%% file: abstract.tex
\begin{abstract}
	We trace the beginning of symbolic dynamics---the study of the shift dynamical system---as it 
	arose from the use of coding to study recurrence and transitivity of geodesics.  It is our assertion
	that neither Hadamard's 1898 paper, nor the Morse-Hedlund papers of 1938 and 1940, which 
	are normally cited as the first instances of symbolic dynamics, truly present the abstract point of
	view associated with the subject today.  Based in part on the evidence of a 1941 letter from 
	Hedlund to Morse, we place the beginning of symbolic dynamics in a paper published by
	Hedlund in 1944.
\end{abstract}

%% file: ShiftText.tex
Symbolic dynamics, in the modern view \cite{Lind-Marcus, Kitchens}, is the dynamical study of the shift automorphism on the space of bi-infinite sequences of symbols, or its restriction to closed invariant subsets.
In this note, we attempt to trace the beginnings of this viewpoint.  While various schemes for symbolic coding of geometric and dynamic phenomena have been around at least since Hadamard (or Gauss: see \cite{Katok-Ugarcovici}), and the two papers by Morse and Hedlund entitled ``Symbolic dynamics'' \cite{Morse-Hedlund:SD1, Morse-Hedlund:SD2} are often cited as the beginnings of the subject, it is our view that the specific, abstract version of symbolic dynamics familiar to us today really began with a paper,``Sturmian minimal sets'' \cite{Hedlund:Sturmian}, published by Hedlund a few years later.  The outlines of the story are familiar, and involve the study of geodesic flows on surfaces, specifically their recurrence and transitivity properties;  this note takes as its focus a letter from Hedlund to Morse, written between their joint papers and Hedlund's, in which his intention to turn the subject into a part of topology is explicit.\footnote{In the interest of full disclosure, we should indicate our possible bias: Hedlund was the first author's dissertation director and hired the second author in his first job.} This letter is reproduced on page \pageref{letter}.

Our focus here is rather narrow, even when it comes to coding geodesic flows.  A recent survey by Katok and Ugarcovici \cite{Katok-Ugarcovici} distinguishes two approaches to such coding: a geometric method, which is our subject, and a second, going back to Gauss and associated with Artin, Koebe and Nielsen which can be regarded as more arithmetic in nature.  
We do not propose to consider the latter in detail.  The survey \cite{Katok-Ugarcovici} discusses technical details of both approaches as well as their subsequent development in recent years.

The beginnings of symbolic dynamics are often traced back to Hadamard's 1898 study of geodesics on surfaces of negative curvature \cite{Hadamard:neg}.  In Part II of this work (\S20), Hadamard gives a coding for the (free homotopy classes of) closed geodesics, essentially as words (up to cyclic permutation) in generators of the fundamental group, and in Part III (\S37) he shows that each word corresponds to a unique closed geodesic.  He then goes on to study unbounded geodesics and in Part VI (\S56) shows that the initial conditions at a point which determine geodesics staying in a bounded region form a perfect, nowhere dense closed set (which is the closure of the conditions yielding closed geodesics).

There are several respects in which Hadamard's paper does not really qualify as a beginning for symbolic dynamics.  First, his coding is limited to \emph{finite} words, coding \emph{closed} geodesics;  he does not appear to envision a coding system encompassing other geodesics.  The bounded, non-periodic geodesics he produces in Part VI are, in passing, seen as determined by a sequence of closed geodesics, but this is not explicitly related to their coding.
Furthermore, Part II of the paper (where the coding is formulated), entitled ``Consid\'erations d'Analysis situs'', is presented as follows:
\begin{quote}
	Ayant reconnu, dans les num\'eros pr\'ec\'edents, l'existence de surfaces \`a courbures oppos\'ees
	et \`a connexion quelconque, nous avons \`a rappeler les principes qui gouvernent l'\'etude 
	des lignes trac\'ees sur de telles surfaces:  principes pos\'es par M. Jordan dans un M\'emoire
	bien connu (2).
\end{quote}
Footnote (2) (``Ce journal, an\'ee 1866'') refers to the second of two back-to-back papers published 33 years before Hadamard in the \textit{Journal de Math{\'e}matiques Pures et Appliqu\'es}
by Jordan \cite{Jordan:Deformation, Jordan:Contours}: the first concerns the role of fundamental contours in determining the homeomorphism type of a surface, and the second presents the notion of a ``class'' of contours (\ie{} free homotopy class) subsequently used by Hadamard in \cite{Hadamard:neg}.  Jordan's notation is that adopted by Hadamard, and he hints at the representation of curves by words (with positive or negative exponents).
Finally, it should be noted that in Hadamard's study, the point of view is geometric rather than dynamic:  geodesics are regarded as oriented curves, and there appears no explicit sense of a ``geodesic flow'';  in particular Hadamard's symbolic coding is static in nature.

In an important paper \cite{Birkhoff:Quelques} published in 1912 (and based on a presentation to the American Mathematical Society in 1909), G. D. Birkhoff analyzes the behavior of recurrent trajectories in a dynamical system defined by a system of ordinary differential equations.  The word ``recurrent'' here corresponds to what we now call ``minimal''.\footnote{What we now call recurrent is called \emph{stable in the sense of Poisson}.}  A collection $M$ of trajectories of a dynamical system is \emph{minimal} if every element of $M$ has all elements of $M$ in its \alp-and \omeg-limit sets;  Birkhoff calls any trajectory belonging to a minimal set ``recurrent''.  He proves that recurrence (in this sense) is equivalent to (what we now call) \emph{almost-periodicity}: for any \epsgo{} there exists a length $T$ such that the whole trajectory is contained in an \eps-neighborhood of any segment of length $T$.  Obvious examples of minimal sets are equilibria and closed orbits;  Birkhoff also notes the example of dense lines on a torus, and calls a recurrent motion \emph{continuous} if the corresponding minimal set forms a continuum of some dimension.  The \apriori{} possibility of \emph{discontinuous} recurrent trajectories is illustrated by the suspension of a nontransitive homeomorphism of the circle with irrational rotation number, and Birkhoff asks whether discontinuous recurrent trajectories can occur in analytical dynamical systems.

Morse, in his 1917 dissertation under Birkhoff (published as \cite{Morse:One-to-One} and \cite{Morse:Recurrent}) establishes the existence of recurrent geodesics of discontinuous type on surfaces of negative curvature and negative Euler characteristic.  He considers the bounded region $S$ obtained by cutting off any infinite ``funnels'' using closed geodesics, and codes the geodesics entirely contained in $S$ by recording the order in which they cross a family of transversals (``normal segments'') that cut $S$ into a simply-connected region---in effect lifting the geodesic to the hyperbolic plane.  He then shows that this coding distinguishes geodesics in $S$, and by constructing the ``Morse sequence'' (discovered earlier and independently by Thue \cite{Thue:Gegenseitige}) proves the existence of discontinuous recurrent geodesics.  Furthermore, he shows by symbolic methods that every closed geodesic in $S$ is a limit of discontinuous recurrent ones.  Since Hadamard had shown that every geodesic in $S$ is a limit of closed geodesics, it follows that the recurrent geodesics of discontinuous type are dense in the set of all geodesics contained in $S$.  Despite the closer connection with dynamical ideas, the point of view in these papers remains geometric:  geodesics are still regarded as curves rather than trajectories, and the coding is used to establish that a geodesic is recurrent (in his sense) and not closed.

In 1920, Birkhoff published a study \cite{Birkhoff:Surface} setting forth a number of ways that the behavior of a dynamical system with two degrees of freedom can be studied by means of the successive intersections of orbits with a transverse surface;  strictly speaking, such a \emph{surface of section} is not entirely transverse to the flow, as it is bounded by closed orbits, but its interior \emph{is} transverse to the flow.  The general setup had been formulated by Poincar\'e in \cite[vol. III, Chap. 27]{Poincare:Nouvelles}, as a means of studying periodic and homoclinic orbits in celestial mechanics. Birkhoff had used the same setup in a limited way for similar purposes in 1917 (\cite{Birkhoff:Two}).
In Chapters 5 and 6 of \cite{Birkhoff:Surface}, Birkhoff goes beyond the study of periodic orbits to study, in some abstraction, the behavior of the ``Poincar\'e map'', defining \alp- and \omeg-limit sets, minimal sets and his related notion of recurrence, as well as transitivity.

In 1924, Artin published a brief but influential paper \cite{Artin:mechanisches} in which he shows that the orbit space of the group of linear fractional transformations with integer coefficients acting on  the hyperbolic plane (in the half-plane model) has a dense geodesic (in fact, the set of these has full measure).  His proof involves a coding of geodesics via the continued fraction expansion of their ``endpoints at infinity'' on the real line.  (His term for transitivity is ``quasiergodicity''.)

In 1927, in the first of his papers on mapping classes \cite{Nielsen:Untersuchungen}, Nielsen formulates a similar coding geometrically, in terms of the fundamental group, to study the axes of hyperbolic transformations on surfaces obtained as quotients of the hyperbolic disc by a Fuchsian group.  Nielsen's approach has some similarities to Morse's coding of geodesics via transverse segments, but the dynamics that comes in is that of the Fuchsian group acting on the universal covering.

In his 1927 book, \emph{Dynamical Systems} \cite{Birkhoff:DynSyst}, a broad survey of work on dynamical systems (primarily of mechanical origin), Birkhoff included Chapter 7, ``General Theory of Dynamical Systems'', which sets forth the notions of wandering and non-wandering orbits, central motions, minimal sets, and transitivity in the general context of the flow generated by a system of differential equations.  Much of this reflected ideas formulated earlier in his 1912 paper \cite{Birkhoff:Quelques}.

In 1935, Birkhoff summed up his work on dynamics in a long paper \cite{Birkhoff:Nouvelles}, ``crowned'' and published by the Pontifical Academy of Sciences.  Chapter 3, a study of behavior near a hyperbolic periodic orbit, is based on a detailed examination of the dynamics of a Poincar\'e map for a transverse section.  By symbolic methods that, several decades later, were modified and used by Smale to prove the ``Smale-Birkhoff'' theorem and to construct the ``horseshoe'', Birkhoff demonstrates the existence of highly complicated ``first return'' behavior for periodic orbits near any orbit homoclinic to a hyperbolic periodic orbit or, more generally, belonging to a loop of heteroclinic connections.

This work forms the background to two papers by Morse and Hedlund, entitled ``Symbolic Dynamics'' \cite{Morse-Hedlund:SD1, Morse-Hedlund:SD2}, published in 1938 and 1940, respectively.  

Hedlund, in his dissertation written under Morse in 1929 \cite{Hedlund:Periodic, Hedlund:Poincare} had proved the existence of a length-minimizing closed geodesic in each free homotopy class for any Riemmannian metric on the torus;  Morse \cite{Morse:Fundamental} had proved the same result for surfaces of higher genus. Hedlund went on to study geodesics on surfaces \cite{Hedlund:Metrical, Hedlund:Metrically, Hedlund:Two, Hedlund:Horocycle}, in particular proving the ergodicity of the geodesic flow on a closed surface of constant negative curvature \cite{Hedlund:Metrical}, using Nielsen's symbolic coding  \cite{Nielsen:Untersuchungen}, and transitivity of the horocycle flow \cite{Hedlund:Horocycle}.
In 1939 he published a survey of results on the dynamics of geodesic flows \cite{Hedlund:Geodesic}, in which he formulates seven types of transitivity, elaborating on Birkhoff's definitions in \cite{Birkhoff:DynSyst}:  these include our notion of topological transitivity (``regional transitivity'', which he notes is equivalent to the existence of a dense trajectory), topological mixing (``permanent regional transitivity''), ergodicity (``metric transitivity'') and mixing (``mixture'') as well as hybrids of topological and ergodic notions of transitivity. He quotes theorems establishing many of these properties for geodesic flows on surfaces of constant negative curvature, as well as an example of a topologically mixing but non-ergodic geodesic flow.   At the end of the article he plugs the work he had started with Morse in \cite{Morse-Hedlund:SD1}:
\begin{quote}
	The development of a symbolic theory apart from its dynamical significance has recently been begun by Morse and the author (cf. Morse [4]).  This initial work includes an extensive analysis of transitive symbolic trajectories.  The full scope of these symbolic methods in dynamics is yet to be determined.
\end{quote}

The first of the Morse-Hedlund papers \cite{Morse-Hedlund:SD1} sets forth a general theory of what we now call shift spaces, focusing on recurrence and transitivity properties of sequences.  The motivation in the introduction refers primarily to geodesic flows on surfaces of negative curvature, but after that the treatment is quite abstract.  The authors' view of the place of their study in dynamics as a whole is stated as follows \cite[pp. 816-817]{Morse-Hedlund:SD1}:
\begin{quote}
	Symbolic dynamics as the authors conceive it forms one of the three divisions
	\begin{enumerate}
		\item representation theory,
		\item symbolic dynamics,
		\item existence of space forms,
	\end{enumerate}
	of the whole theory.  The representation theory is concerned with the conditions on space forms 
	under which trajectories admit a one-to-one symbolic representation in terms of which the 
	recurrence or transitivity of the trajectory can be determined.  These conditions will involve
	the Poincar\'e fundamental group of the space and differential conditions such as that of
	uniform instability (cf. Morse [4]\footnote{This reference is our \cite{Morse:One-to-One}.}, p. 64).
	In (3) one is concerned with the existence of space forms satisfying the conditions discovered
	in (1).  The questions involved are rather deep extensions of the Hilbert, Koebe theory of spaces
	of negative curvature (cf. Hilbert [1]\footnote{\cite{Hilbert:Flachen}}, 
	and Koebe [1]\footnote{\cite{Koebe:Riemannsche}}).  A simple typical theorem is that there exists
	no two-dimensional Riemannian manifold of the topological type of the torus satisfying the
	condition of uniform geodesic instability.  The bearing of such studies on questions of topological
	and metric transitivity will be made clear in later papers.
\end{quote} 
Clearly, Morse and Hedlund view their paper as initiating a new branch of the theory of dynamical systems.  However, it does not seem to us that the \emph{shift dynamical system} is as yet considered as an object of study.

Beginning with a finite alphabet, Morse and Hedlund define an \emph{$I$-trajectory}  to be a 
two-sided indexed sequence of letters; a \emph{symbolic element} $E(r,a)$ is an $I$-trajectory $a=...\as{-1}\as{0}\as{1}...$ together with a choice of a distinguished position $r$ on it.  The space of all symbolic elements is given the metric
\begin{equation*}
	d(E(r,a),E(s,b)=\recip{m}
\end{equation*}
when $\as{r-m}...\as{r+m}$ and $\bs{s-m}...\bs{s+m}$ are the longest symmetric words centered on the distinguished positions which agree termwise (elements whose distinguished positions have different values are at infinite distance).  They establish that this gives the space of symbolic elements the topology of a Cantor set.  The space of $I$-trajectories is given the metric
\begin{equation*}
	[a,b]=\limsup_{n\to\infty}\recip{2n+1}\sum_{-n}^{n}\delta(\as{i},\bs{i})
\end{equation*}
(where $\delta(\as{i},\bs{i})=1$ or $0$, as \as{i} and \bs{i} are the same or different) which they view as an analogue of the sup metric on functions, as used by Besicovitch in his treatment of almost-periodic functions \cite{Besicovitch}---in fact, their notation closely follows his.  They define an $I$-trajectory $a$ to be \emph{almost periodic} if for every \epsgo{} the iterates $D_{r}$ of the shift automorphism which satisfy
\begin{equation*}
	[D_{r}(a),a]<\eps
\end{equation*}
form a relatively dense set of integers (that is, there is an integer $N$ such that any set of $N$ consecutive integers intersects the set).
Note that this is stronger than what we now call almost-periodicity, as there is a uniformity condition involved.  They consider  subsets of the space of trajectories defined by admissible blocks;  their admissibility rules appear to be of finite type, although they state a family of conditions \cite[p. 823]{Morse-Hedlund:SD1} which are far more restrictive, and appear to be motivated by Nielsen's formalism, to which they explicitly refer as an example.  Again, they show that the subspace so defined has the topology of a Cantor set. They then study limit trajectories and minimal sets of trajectories from a symbolic point of view, and present the Morse sequence.  The last 60\% of the paper (pp. 833-864) is taken up with a number of functions that measure the ``speed'' of recurrence;  these need not occupy us in detail here.  It should be noted that, despite the dynamical background, and the appearance of the shift automorphism in two places (pp. 817 and 822), it is used in a way analogous to Besicovitch's use of translations to study almost-periodic functions (in fact, as we have noted, their notation is the same);   there is no sense of a dynamical system generated by iteration of the shift.

Morse and Hedlund's second paper \cite{Morse-Hedlund:SD2} concerns a specific class of subshifts, which they explain characterize the geodesics on a flat torus.  These are built on an alphabet of two symbols and are defined by the condition that for each symbol, any two maximal blocks of consecutive appearances of the symbol differ in length by at most one.  At the end of their previous paper, they had noted the relation of this condition to the Sturm Separation Theorem concerning the distribution of zeroes of the solution of the differential equation
\begin{equation*}
	\ypp+f(x)y=0
\end{equation*}
where \fof{x} is periodic with period one: one symbol represents the locations of zeroes, the other the locations of integers (it is assumed without loss of generality that the solution has no integer zeroes).
They call such trajectories \emph{Sturmian trajectories}.  This paper is a detailed algebraic study of various combinatorial functions that characterize a Sturmian trajectory.  Again, there is no explicit dynamical system here.

The journal lists \cite{Morse-Hedlund:SD2} as received June 19, 1939.  Two years later, Hedlund wrote to Morse as follows\label{letter}:
\begin{quote}
	\flushright{Charlottesville, Virginia}\\
	\flushright{June 7, 1941}
	
	\flushleft{Dear Marston:}\\
\smallskip	
\raggedright{	As you probably know, a number of topologists are becoming interested in a study of the structure of the orbits obtained when a topological transformation is iterated on a space $X$, say a separable metric space.  They term an orbit the set $T^{n}(x)$, $n=...-1,0,+1,+2,...$, where $x$ is a point of $X$, define a point $x$ to be periodic if there exists an integer $m>0$ such that $t^{m}(x)=x$, and define apoint $x$ to be almost periodic if there exists a sequence of integers $\ns{1}<\ns{2}<...$ such that $\lim_{i\to\infty}T^{\ns{i}}(x)=x$.  This last definition is of course the well known property of (positive) stability in the sense of Poisson and the term almost periodic is somewhat of a misnomer.  These topologists are not in the least aware that there is an immmense amount of material in dynamics which they should know and they will probably rediscover such interesting things as minimal sets, recurrent motions, minimal centers of attraction, central motions, transitivity, permanent regional transitivity, etc., in the not too distant future.  For example, in the last issue of Mathematical Reviews (see page 179, review of a paper of Schweigert) Ayres comments that Schweigert has an interesting example for which the periodic points are everywhere dense in a space, but not all the points of the space are periodic.  Now this is such a common occurrence in dynamics that we scarcely wonder at it any more.\\
\smallskip	
	But I wonder if it is their fault that these things are not better known.  For a person who hasn't dealt considerably with these matters it might be a hard task to dig the material out of the literature.  One reason is perhaps that in dynamics we deal largely with flows, whereas the topologists deal with the discrete case of a single transformation and its iterates.  Though the two are not essentially different, it seems to me that something should be done about this situation.  Yet I hesitate about publishing material which can onlybe [\textit{sic}] considered a rehash of mathematics which is well known ( to at least a dozen people ).\\
	\smallskip
	However, one simple example occurs to me which might make for more awareness concerning the results of dynamics.  Let $M$ be the space of symbolic elements of our first paper on symbolic dynamics and let $T$ be the transformation which shifts the index by one, say to the right.  Here is a topological transformation which ought to be complex enough to suit the heart of even the most pathological topologist.  The periodic trajectories are everywhere dense;  the non-periodic recurrent trajectories are everywhere dense;  there are transitive trajectories and they form a residual set;  the non-periodic, non-recurrent, non-transitive trajectories which are stable in the sense of Poisson are everywhere dense;  there are trajectories asymptotic to almost anything;  the transformation is permanently regionally transitive.  What do you think of giving them this example on which to chew?  In view of what we have avaible[\textit{sic}]  in SDI, it should not occupy much space.\\
	\smallskip
	Though the preceding example is a good one, it has one defect.  The space $M$ is, as we showed, compact, perfect, and totally disconnected.  The last property of being totally disconnected should not be essential to the situation and is not characteristic of classical dynamical systems, where the underlying spaces are manifolds.  The space $M$ is disconnected because the metric which we chose to topologize the space assumes only discrete values.  Would it be possible to topologize the space in some other fashion so that it becomes say a continuum (compact, connected space)?  I began thinking of this last night and the answer hasn't occurred to me yet.  It should be possible to define a non-trivial metric in the space $M$ of elements such that it goes to zero as larger and larger blocks with center at the preferred symbol become identical and yet such that $M$ is connected.  It may be necessary to identify a denumerable set of elements in pairs, but that would not be objectionable.  I should think that almost any sort of a space might be obtained with the proper choice of metric, and the whole problem seems interesting.\\
	\smallskip
	If all this is of interest to you, I will be glad to learn your reactions.
}\\
\smallskip
\flushright{As ever,}\\
\flushright{G A H}
\end{quote}

Several features of this letter deserve comment.  The first thing that strikes one is the penultimate paragraph, in which Hedlund wonders whether by a different choice of metric the space of symbolic elements could be made connected.  Of course, on closer reading this is modified by the comment that this might be accomplished by  identifying a countable number of pairs of points.  Anachronistically, one could view this as a precursor of the construction of Markov partitions \cite{Parry:Symbolic, Berg, Adler-Weiss1, Adler-Weiss2}.  Of course, the fact that the shift space itself has the topology of a Cantor set is these days taken for granted, and not viewed as problematic: the coding is a map from an initially given dynamical system to the shift dynamical system.  Second, Hedlund is fully aware of the idea of a discrete dynamical system, and sees it as fundamentally equivalent to the idea of a flow.  He is also completely aware of the interplay of topological and dynamical features of the system, something that is not clear in the earlier joint papers.  Finally, the implied distinction between dynamics and topology has been somewhat erased in more recent years.

Perhaps a comment on the interests of ``a number of topologists'' mentioned in the first sentence of Hedlund's letter is in order.  A search of the \emph{Mathematical Reviews} (which began only a year before Hedlund's letter, in 1940) and \emph{Zentralblatt} (begun about ten years earlier) of that period reveals a large number of works concerning periodic and fixedpoint behavior of iterated transformations on manifolds and metric spaces.  Of course, there were some more sophisticated precursors:  Brouwer's characterization of fixedpoint-free transformations of the plane (the ``Brouwer translation theorem'' \cite{Brouwer:Translation}), Denjoy's work on flows on the torus or, equivalently, diffeomorphisms of the circle \cite{Denjoy:Courbes, Denjoy:Caracteristiques}, (see also \cite{Kampen:Topological}), and the work of Fatou \cite{Fatou1, Fatou2, Fatou3} and Julia \cite{Julia} on iteration of rational functions.  

Hedlund's reference to ``the last issue of Mathematical Reviews'' may give the impression that his acquaintance with the work of ``some topologists'' was second-hand.  The paper trail suggests otherwise.
The review in question, a one-paragraph review of \cite{Schweigert} by  Ayres: \textit{Math. Rev.} 2, 179b (MR 3198), begins by referring to a paper by Hall and  Schweigert \cite{Hall-Schweigert} which is being generalized by Schweigert in the paper under review.  Hall was no stranger to Hedlund.  A paper by Hall and Kelley \cite{Hall-Kelley} which appears in Hedlund's bibliography
to \cite{Hedlund:Sturmian}, was published in 1941, based on a presentation to the American Mathematical Society in September 1939.  This paper concerns variants of periodicity for an iterated self-homeomorphism of a compact metric space, including uniform and non-uniform versions of almost-periodicity.  Hall and Kelley formulate the notion of a minimal set (which they call ``irreducibly fixed'') and show that this is equivalent to every orbit being dense, and that either such a  set is a single periodic orbit, or every orbit in it is almost periodic.  In a footnote, they acknowledge that ``It has been pointed out to the authors that...[these results]...are precisely analogous to certain results of G. D. Birkhoff for continuous flows...''  In an earlier footnote, they achnowledge
\begin{quote}
	This paper was started when the authors were in residence at the University of Virginia, 
	the first named author as a National Research Fellow.
\end{quote}
Recall that Hedlund's letter was written in Charlottesville, Virginia, where he had joined the University of Virginia faculty in 1939.

In any case, it appears that Morse failed to respond to Hedlund's letter in any substantial way.  Hedlund wrote a new article, ``Sturmian Minimal Sets'' \cite{Hedlund:Sturmian}, submitted to the journal in January, 1944.  The minimal sets of the title are addressed in the second half of the paper.  Hedlund begins boldly, explicitly introducing the notion of an orbit (note: not ``trajectory'') and semi-orbit for a discrete dynamical system.  He formulates discrete dynamical system versions of the definitions in \cite{Birkhoff:DynSyst} of $\alpha$- and $\omega$-limit sets, invariant sets, minimal sets, recurrent orbit (in the sense of being contained in a minimal set), and almost-periodic orbit, and notes the equivalence of the last two notions.  Then he comments that the dual terminology ``recurrent'' and ``almost periodic'' for equivalent notions is redundant;  he argues that the latter is the better terminology for this notion, and suggests that ``recurrent'' be saved for ``Poisson stable'' (our current meaning of ``recurrent'').
He then repeats the definitions of symbolic trajectory and symbolic element from \cite{Morse-Hedlund:SD1}, and introduces a modified version of the metric on symbolic elements (replacing $m$ with $m+1$, so that elements which agree only at the distinguished position are at distance $1$), noting that this gives the space of symbolic elements the topology of a Cantor set---apparently abandoning his concerns in the letter about this.

Then, significantly, Hedlund introduces the shift map $S$ and proves that it is a one-to-one, continuous transformation of the space of symbolic elements to itself.  He notes that the symbolic trajectories are in one-to-one correspondence with the orbits of $S$, with periodic trajectories corresponding to periodic orbits, and recalls the existence of the Morse sequence and Sturmian minimal trajectories of \cite{Morse-Hedlund:SD2} as examples of almost-periodic non-periodic orbits.   He then proceeds to construct Sturmian minimal sets by coding orbits of a  rotation of the circle by $\beta$ radians, where $\beta$ is an irrational multiple of $\pi$, using a partition into arcs of length $\beta$ and $2\pi-\beta$, respectively.\footnote{Hedlund here works with translation on the real line, but his coding depends only on the mod 1 positions of points.}  He associates to each orbit two symbolic elements, corresponding to making the atoms of the partition right- or left-open.  He proceeds to prove first that the symbolic sequences which arise this way correspond to almost-periodic orbits of the rotation, and then that this corresponds to being minimal under $S$.  He proves that the minimal set (for $S$) so obtained is compact, perfect and totally disconnected, and contains a pair of doubly asymptotic orbits.  He then defines the notion of an orbit-preserving transformation with respect to a given discrete dynamical system (this is our notion of a self-conjugacy) and proves that in his minimal sets it is not always possible to find an orbit-preserving transformation taking one orbit in the set to an arbitrarily designated second orbit in the set.  He attributes the corresponding question for flows to Birkhoff.   He defines a notion of almost-periodicity for a transformation (as opposed to a single orbit) and shows that the restriction of $S$ to his minimal sets does not have this property.  However, he then defines a notion of \emph{local} almost-periodicity, and shows that his minimal sets do have this property.  Finally, he defines a minimal set to be \emph{powerfully minimal} if it is minimal under all nonzero iterates of the discrete dynamical system, and proves that his minimal set has this property.

While the focus in this paper is on a specific class of minimal sets (the ``Sturmian'' ones), Hedlund's letter suggests that these are now a case study for a more  general, abstract study of minimal sets and dynamical properties of the discrete dynamical system defined by $S$ on the space of symbolic elements, and thus open the door to the branch of topological dynamics we now call ``symbolic dynamics''.